\documentclass[leqno]{amsart}
\usepackage{amsthm}
\usepackage{color}
\usepackage{amssymb}
\usepackage{amsmath,amsfonts,enumerate}
\usepackage{amsthm}

\numberwithin{equation}{section}

\newtheorem{thm}{Theorem}
\newtheorem{lem}[thm]{Lemma}

\newtheorem{definition}[thm]{Definition}
\newtheorem{corollary}[thm]{Corollary}

\newtheorem{proposition}[thm]{Proposition}

\begin{document}

\title{The non-commutative Hardy-Littlewood maximal function on symmetric spaces of $\tau-$measurable operators
}

\author{Y. Nessipbayev}
\address{
c;
Al-Farabi Kazakh National University, 050040 Almaty, Kazakhstan;
International Information Technology University, 050040 Almaty, Kazakhstan;}
\email{yerlan.nessipbayev@gmail.com}

\author{K. Tulenov}
\address{
Al-Farabi Kazakh National University, 050040 Almaty, Kazakhstan;
Institute of Mathematics and Mathematical Modeling, 050010 Almaty, Kazakhstan.}
\email{tulenov@math.kz}

\date{}
\begin{abstract} In this paper, we investigate the Hardy-Littlewood maximal function on non-commutative symmetric spaces. We complete the results of T. Bekjan and J. Shao. Moreover, we refine the main results of the papers \cite{Bek} and \cite{Sh}.
\end{abstract}

\maketitle

\section{Introduction}
The Hardy-Littlewood maximal operator (or function) $M$ is an important non-linear operator with numerous applications in real analysis and harmonic analysis. It takes a locally integrable function $f : \mathbb{R}^d \to \mathbb{C}$ and returns another function $Mf$ that, at each point $x \in \mathbb{R}^d$, gives the maximum average value that $f$ can have on balls centered at that point. More precisely,

\begin{equation}\label{MF}
Mf(x)=\sup_{r>0}\frac{1}{|B(x,r)|}\int_{B(x,r)}|f(y)|dy,
\end{equation}
where $B(x,r)$ is the ball of radius $r$ centred at $x$, and $|E|$ denotes the $d$-dimensional Lebesgue measure of $E \subset \mathbb{R}^d.$
There is an uncountable number of papers devoted to investigation of the Hardy-Littlewood function defined by the formula \eqref{MF}. For instance, see \cite{BSh, Gra, SW} and references therein.

It is well known that there are positive constants $c_1$ and $c_2$, depending only on $d$, such that

\begin{equation}\label{maximal cesaro}
c_1 \mu(Mf) \leq C \mu(f) \leq c_2 \mu(Mf)
\end{equation}
for every locally integrable function $f$ on $\mathbb{R}^d$ (see \cite[Theorem III. 3.8, p.122]{BSh}), where $C$ is the Ces\`{a}ro operator defined by \eqref{cesaro} below and $\mu(f)$ is the decreasing rearrangement of the function $|f|.$

In this paper we mainly deal with the non-commutative version of this result. In fact, for our purposes it is sufficient to consider only the first inequality in \eqref{maximal cesaro}.
To be precise, in this work we study the non-commutative Hardy-Littlewood maximal function on symmetric spaces of $\tau-$measurable operators.
Non-commutative maximal inequalities were studied, in particular, in works of M. Junge and Q. Xu \cite{Junge, Xu}. Later T. Mei \cite{Mei} presented a version of the non-commutative Hardy-Littlewood maximal inequality for an operator-valued function. Another version of the (non-commutative) Hardy-Littlewood maximal function was introduced by T. Bekjan \cite{Bek}. Also, in his work T. Bekjan defined the (non-commutative) Hardy-Littlewood maximal function for $\tau-$measurable operators and, among other things, obtained weak $(1,1)-$type and $(p,p)-$type inequalities for the Hardy-Littlewood maximal function. Later J. Shao investigated the Hardy-Littlewood maximal function on non-commutative Lorentz spaces associated with finite atomless von Neumann algebra (see \cite{Sh}). Using the techniques in \cite{Gra}, J. Shao obtained the $(p,q)-(p,q)-$type inequality for the Hardy-Littlewood maximal function on non-commutative Lorentz spaces. Namely, for an operator $T$ affiliated with a semi-finite von Neumann algebra $\mathcal{M},$ the Hardy-Littlewood maximal function of $T$ is defined by
$$MA(x)=\sup_{r>0}\frac{1}{\tau(E_{[x-r,x+r]}(|A|))}\tau(|A|E_{[x-r,x+r]}(|A|)), \,\ x\geq 0.$$
While the classical Hardy-Littlewood maximal function of a Lebesgue measurable function $f: \mathbb{R} \to \mathbb{R}$, denoted by $Mf(x)$, is defined as
$$Mf(x)=\sup_{r>0}\frac{1}{m([x-r,x+r])}\int_{[x-r,x+r]}|f(t)|dt,$$
where $m$ is a Lebesgue measure on $(-\infty, \infty)$ \cite{SW}.
In view of spectral theory, $|A|$ is represented as
$$|A|=\int_{\sigma(|A|)} tdE_t,$$
and $MA(|A|)$ is represented as $MA(x)$. Thus, for the operator $A,$ Bekjan's consideration is that $MA(|A|)$ is defined as the operator analogue of the Hardy-Littlewood maximal function in the classical case. Our purpose is to investigate the non-commutative Hardy-Littlewood maximal function $M$ in the sense of T. Bekjan (see \cite{Bek}).

In this paper we obtain upper estimate of generalized singular number of the non-commutative Hardy-Littlewood maximal function (see Section 3). In particular, we show that the generalized singular number of the Hardy-Littlewood maximal function is estimated from above by the Ces\`{a}ro operator, which we address in Theorem \ref{main th} below.

We extend the results obtained by J. Shao (see \cite{Sh}) from a finite von Neumann algebra $\mathcal{M}$ to a general semi-finite $\mathcal{M}.$ Also, while J. Shao considered the non-commutative Hardy-Littlewood function on Lorentz spaces, we rather show its boundedness from one non-commutative symmetric space to another (see Theorem \ref{main1} below), including its boundedness from a non-commutative Lorentz space to another (see Proposition \ref{16} below) and from a non-commutative Marcinkiewicz space to another (see Proposition \ref{17} below). We show it using an approach similar to the commutative case. The key results, which we apply, may be found in \cite{Bek}. Finally, we provide interesting examples showing that the non-commutative Hardy-littlewood maximal function is bounded from one particular non-commutative Lorentz space to another, and from one particular non-commutative Marcinkiewicz space to another as well (see Section 4).

\section{Preliminaries}
Let $(\mathbb{R_+},m)$ denote the measure space $\mathbb{R_+} = (0,\infty)$ equipped with Lebesgue measure $m.$
 Let $L(\mathbb{R_+},m)$ be the space of all measurable real-valued functions on $\mathbb{R_+}$ equipped with Lebesgue measure $m$ i.e. functions which coincide almost everywhere are considered identical. Define $S(\mathbb{R_+},m)$ to be the subset of $L(\mathbb{R_+},m)$ which consists of all functions $x$ such that $m(\{t : |x(t)| > s\})$ is finite for some $s > 0.$
 For $x\in S(\mathbb{R_+})$ we denote by $\mu(x)$ the decreasing rearrangement of the function $|x|.$ That is,
$$\mu(t,x)=\inf\{s\geq0:\ m(\{|x|>s\})\leq t\},\quad t>0.$$

We say that $y$ is submajorized by $x$ in the sense of Hardy--Littlewood--P\'{o}lya (written $y\prec\prec x$) if
$$\int_0^t\mu(s,y)ds\leq\int_0^t\mu(s,x)ds,\quad t\geq0.$$

Let $\mathcal{M}$ be a semifinite von Neumann algebra on a separable Hilbert space $H$ equipped with a faithful normal semifinite trace $\tau.$
A closed and densely defined operator $A$ affiliated with $\mathcal{M}$ is called $\tau$-measurable if $\tau(E_{|A|}(s,\infty))<\infty$ for sufficiently large $s.$ We denote the set of all $\tau$-measurable operators by
$S(\mathcal{M},\tau).$ Let $Proj(\mathcal{M})$ denote the lattice of all projections in $\mathcal{M}.$ For every $A\in S(\mathcal{M},\tau),$ we define its singular value function $\mu(A)$ by setting
$$\mu(t,A)=\inf\{\|A(1-P)\|_{\mathcal{L}_{\infty}(\mathcal{M})}:P\in Proj(\mathcal{M}),\quad \tau(P)\leq t\}, \quad t>0,$$
where the norm $\|\cdot\|_{\mathcal{L}_{\infty}(\mathcal{M})}$ is the usual operator (uniform) norm.
Equivalently, for positive self-adjoint operators $A\in S(\mathcal{M},\tau),$ we have
$$n_A(s)=\tau(E_A(s,\infty)),\quad \mu(t,A)=\inf\{s: n_A(s)<t\}, \quad t>0.$$
 An operator in $S(\mathcal{M},\tau)$ is called $\tau$-compact if $\mu(\infty,A)=0.$ This notion is a direct generalization of the ideal of compact operators on a Hilbert
space $H.$ For more details on generalised singular value functions and $\tau$-compact operators, we refer the reader to \cite{FK} and \cite{LSZ}.
Let $\mathcal{L}_{loc}(\mathcal{M},\tau)$ be the set of all $\tau$-measurable operators such that
\begin{equation}\label{*}\tau(|A|E_{I}(|A|))<+\infty,
\end{equation}
for all bounded intervals $I\subset[0,+\infty)$.

We have for $A,B \in S(\mathcal{M},\tau)$ (see for instance \cite[Corollary 2.3.16 (a)]{LSZ})
\begin{equation}\label{triangle svf}
\mu(t+s,A+B)\leq\mu(t,A)+\mu(s,B),\quad t,s>0.
\end{equation}

If $A,B\in S(\mathcal{M},\tau),$ then we say that $B$ is submajorized by $A$ (in the sense of Hardy--Littlewood--P\'{o}lya) denoted by $\mu(B)\prec\prec\mu(A),$ that is,
$$\int_0^t\mu(s,B)ds\leq\int_0^t\mu(s,A)ds,\quad t\geq0.$$


\begin{definition} \cite[Definition 1]{Bek} For $A\in \mathcal{L}_{loc}(\mathcal{M},\tau)$, we define the maximal function of $A$ by
	$$MA(x)=\sup_{r>0}\frac{1}{\tau(E_{[x-r,x+r]}(|A|))}\tau(|A|E_{[x-r,x+r]}(|A|)), \,\ x\geq 0,$$
	(let $\frac{0}{0}=0)$. $M$ is called the non-commutative Hardy-Littlewood maximal function.
\end{definition}

$MA(|A|)$ is represented as $MA(x).$ Then for $A,$ we consider $MA(|A|)$ as the operator analogue
of the Hardy-Littlewood maximal function in the classical case.  Hence roughly speaking, $MA(|A|)$ stands in relation to $A$ as $Mf(x)$
stands in relation to $f$ in classical analysis, i.e.
\begin{equation}\label{maxop}
MA(|A|)=\int_{\sigma(|A|)}MA(\lambda)dE_{\lambda}(|A|),
\end{equation}
where $\sigma(|A|)$ is the spectrum of $|A|.$

\subsection{Symmetric (Quasi-)Banach Function and Operator Spaces}

\begin{definition}\label{Sym} We say that $(E,\|\cdot\|_E)$ is a symmetric (quasi-)Banach function space on $I$ if the following hold:
\begin{enumerate}[{\rm (a)}]
\item $E$ is a subset of $S(\mathbb{R_+},m);$
\item $(E,\|\cdot\|_E)$ is a (quasi-)Banach space;
\item If $x\in E$ and if $y\in S(\mathbb{R_+},m)$  are such that $|y|\leq|x|,$ then $y\in E$ and $\|y\|_E\leq\|x\|_E;$
\item If $x\in E$ and if $y\in S(\mathbb{R_+},m)$ are such that $\mu(y)=\mu(x),$ then $y\in E$ and $\|y\|_E=\|x\|_E.$
\end{enumerate}
\end{definition}
For the general theory of symmetric (quasi-)Banach function spaces,
we refer the reader to \cite{BSh,KPS,LT}.

We restrict our attention to the subclass of symmetric {\it Banach spaces $E$ with Fatou norm}, that is, when the norm closed unit ball $B_E$ of $E$ is closed in $E$ with respect to almost everywhere convergence.

Let $\mathcal{M}$ be a semifinite von Neumann algebra on a Hilbert space $H$ equipped with a faithful normal semifinite trace $\tau.$ The von Neumann algebra
$\mathcal{M}$ is also denoted by $L_{\infty}(\mathcal{M}),$ and for all $A\in L_{\infty}(\mathcal{M}),$ $\|A\|_{L_{\infty}(\mathcal{M})}:=\|A\|$.
Moreover, $\|A\|_{L_{\infty}(\mathcal{M})}=\|\mu(A)\|_{L_{\infty}(\mathbb{R_+})}=\mu(0,A)$, $A\in L_{\infty}(\mathcal{M})$.

\begin{definition}\label{NC Sym} Let $E(\mathcal{M},\tau)$ be a linear subset in $S(\mathcal{M},\tau)$ equipped with a complete quasi-norm $\|\cdot\|_{E(\mathcal{M},\tau)}.$ 
We say that $E(\mathcal{M},\tau)$ is a symmetric operator space (on $\mathcal{M}$, or in $S(\mathcal{M},\tau)$) if
for $A\in E(\mathcal{M},\tau)$ and for every $B\in S(\mathcal{M},\tau)$ with $\mu(B)\leq\mu(A),$ we have $B\in E(\mathcal{M},\tau)$ and $\|B\|_{E(\mathcal{M},\tau)}\leq\|A\|_{E(\mathcal{M},\tau)}$.

A symmetric function space is the term reserved for a symmetric
operator space when $\mathcal{M}=L_{\infty}(\mathbb{R_+},m)$.
\end{definition}

Recall the construction of a symmetric (quasi-)Banach operator space (or non-commuative symmetric (quasi-)Banach space) $E(\mathcal{M},\tau).$ Let $E$ be a symmetric (quasi-)Banach function (or sequence) space on $(0,\infty)$. Set
$$E(\mathcal{M},\tau)=\Big\{A\in S(\mathcal{M},\tau):\ \mu(A)\in E\Big\}.$$
We equip $E(\mathcal{M},\tau)$ with a natural norm
$$\|A\|_{E(\mathcal{M},\tau)}=\|\mu(A)\|_E,\quad A\in E(\mathcal{M},\tau).$$
This is a (quasi-)Banach space with the (quasi-)norm $\|\cdot\|_{E(\mathcal{M},\tau)}$ and is called the (non-commutative) symmetric operator space associated with $(\mathcal{M},\tau)$ corresponding to $(E,\|\cdot\|_{E})$.
An
extensive discussion of the various properties of such spaces can be found
in \cite{KS,LSZ} (see also \cite{STZ1,STZ2}).
Futhermore, the following fundamental theorem was proved in \cite{KS} (see also \cite[Question 2.5.5, p. 58]{LSZ}).
\begin{thm} Let $E$ be a symmetric function space on $\mathbb{R_+}$ and let $\mathcal{M}$ be a semifinite von Neumann algebra. Set
$$E(\mathcal{M},\tau)=\Big\{A\in S(\mathcal{M},\tau):\ \mu(A)\in E\Big\}.$$
So defined $(E(\mathcal{M},\tau),\|\cdot\|_{E(\mathcal{M},\tau)})$
is a symmetric operator space.
\end{thm}

For simplicity, in the future we will denote $E(\mathcal{M},\tau)$ by $E(\mathcal{M})$.

For $1\leq p<\infty,$ we set
$$L_p(\mathcal{M})=\{A\in S(\mathcal{M},\tau):\ \tau(|A|^p)<\infty\},\quad \|A\|_{L_p(\mathcal{M})}=(\tau(|A|^p))^{\frac1p}.$$
The Banach spaces $(L_p(\mathcal{M}),\|\cdot\|_{L_p(\mathcal{M})})$ ($1\leq p<\infty$) are separable.
Moreover, when $p=2$ this space
becomes the Hilbert space with the inner product
$$<A,B>:=\tau(B^{*}A),\,\, A,B\in L_{2}(\mathcal{M}),$$  where $B^{*}$ is adjoint operator of $B.$
It is easy to see that these spaces are examples of symmetric spaces.

Throughout this paper, we write $\mathcal{A}\lesssim \mathcal{B}$ if there is a constant $c_{abs}> 0$ such that
$\mathcal{A}\leq c_{abs}\mathcal{B}.$ We write $\mathcal{A}\approx \mathcal{B}$ if both $\mathcal{A}\lesssim \mathcal{B}$ and $\mathcal{A}\gtrsim \mathcal{B}$ hold, possibly with different
constants.

\subsection{$L_{1}\cap L_{\infty}$ and $L_{1}+L_{\infty}$ spaces}

Two examples below are of particular interest. Consider the separated topological vector space $S(\mathbb{R_+})$ consisting of all measurable functions $x$ such that $m(\{t : |x(t)| > s\})$ is finite for some $s > 0$ with the topology of convergence in measure. Then the spaces $L_{1}(\mathbb{R_+})$ and $L_{\infty}(\mathbb{R_+})$ are algebraically and topologically imbedded in the topological vector space $S(\mathbb{R_+})$, and so these spaces form a Banach couple (see \cite[Chapter I]{KPS} for more details). The space $(L_{1}\cap L_{\infty})(\mathbb{R_+})=L_{1}(\mathbb{R_+})\cap L_{\infty}(\mathbb{R_+})$ consists of all bounded summable functions $x$ on $\mathbb{R_+}$ with norm
$$\|x\|_{(L_{1}\cap L_{\infty})(\mathbb{R_+})}=\max\{\|x\|_{L_{1}(\mathbb{R_+})},\|x\|_{L_{\infty}(\mathbb{R_+})}\}, x\in(L_{1}\cap L_{\infty})(\mathbb{R_+}).$$
The space $(L_{1}+L_{\infty})(\mathbb{R_+})=L_{1}(\mathbb{R_+})+ L_{\infty}(\mathbb{R_+})$ consists of functions which are sums of bounded measurable and summable functions $x\in S(\mathbb{R_+})$
equipped with the norm given by
$$\|x\|_{(L_{1}+L_{\infty})(\mathbb{R_+})}=\inf\{\|x_1\|_{L_{1}(\mathbb{R_+})}+\|x_2\|_{L_{\infty}(\mathbb{R_+})}:x=x_1+x_2, $$
$$\, x_{1}\in L_{1}(\mathbb{R_+}), \, x_2\in L_{\infty}(\mathbb{R_+})\}.$$
For more details we refer the reader to \cite[Chapter I]{BSh},\cite[Chapter II]{KPS}.
We recall that that every symmetric Banach function space on $\mathbb{R_+}$ (with respect to Lebesgue measure) satisfies
$$(L_1\cap L_{\infty})(\mathbb{R_+})\subset E(\mathbb{R_+})\subset (L_1+L_{\infty})(\mathbb{R_+})$$
with continuous embeddings (see for instance \cite[Theorem II. 4.1. p. 91]{KPS}).

We define the space $L_{1}(\mathcal{M})+L_{\infty}(\mathcal{M})$ as the class of those operators $A\in S(\mathcal{M},\tau)$ for which
\begin{eqnarray*}\begin{split}\|A\|_{L_{1}(\mathcal{M})+L_{\infty}(\mathcal{M})}
&:=\inf\{\|A_1\|_{L_1(\mathcal{M})}+\|A_2\|_{L_{\infty}(\mathcal{M})}:\\
&A=A_1+A_2,\, A_1\in L_{1}(\mathcal{M}), A_2\in L_{\infty}(\mathcal{M})\}<\infty.
\end{split}\end{eqnarray*}
\begin{thm}\cite[Theorem III. 9.16, p. 96]{DdePS}\label{s-calcul thm} If $A\in S(\mathcal{M},\tau),$ then
\begin{equation}\label{L+L}
\|A\|_{L_{1}(\mathcal{M})+L_{\infty}(\mathcal{M})}=\int_{0}^{1}\mu(t,A)dt, \,\,\,\,\, A\in L_{1}(\mathcal{M})+L_{\infty}(\mathcal{M}).
\end{equation}
\end{thm}

A standard argument shows that $L_{1}(\mathcal{M})+L_{\infty}(\mathcal{M})$ is
a Banach space with respect to this norm, which is continuously embedded in $S(\mathcal{M},\tau).$ It follows from Theorem \ref{s-calcul thm} that if $A\in S(\mathcal{M},\tau),$ then $A\in
L_{1}(\mathcal{M})+L_{\infty}(\mathcal{M})$ if and only if $\int_{0}^{t}\mu(s,A)ds<\infty$ for all $t>0.$

In particular, if $A\in S(\mathcal{M},\tau),$ then $A\in L_{1}(\mathcal{M})+L_{\infty}(\mathcal{M})$ if and only if $\mu(A)\in (L_{1}+L_{\infty})(\mathbb{R_+})$ and
$$\|A\|_{L_{1}(\mathcal{M})+L_{\infty}(\mathcal{M})}=\|\mu(A)\|_{(L_{1}+L_{\infty})(\mathbb{R_+})}, \,\,\,\,\, A\in L_{1}(\mathcal{M})+L_{\infty}(\mathcal{M}).$$
The space $L_{1}(\mathcal{M})+L_{\infty}(\mathcal{M})$ is also denoted by $(L_{1}+L_{\infty})(\mathcal{M}).$
Similarly, define the intersection space of Banach spaces $L_{1}(\mathcal{M})$ and $L_{\infty}(\mathcal{M})$ as follows.
$$(L_{1}\cap L_{\infty})(\mathcal{M})=\{A\in S(\mathcal{M},\tau): \,\ \|A\|_{(L_{1}\cap L_{\infty})(\mathcal{M})}<\infty\}, $$
where the norm on $(L_{1}\cap L_{\infty})(\mathcal{M})$ is defined by
$$
\|A\|_{(L_{1}\cap L_{\infty})(\mathcal{M})}=\max\{\|A\|_{L_{1}(\mathcal{M})},\|A\|_{L_{\infty}(\mathcal{M})}\}, \,\,\  A\in (L_{1}\cap L_{\infty})(\mathcal{M}).
$$
It is easy to see that $(L_{1}\cap L_{\infty})(\mathcal{M})$ is a Banach space with respect to this norm. It should be observed that if $A\in S(\mathcal{M},\tau),$ then $A\in (L_{1}\cap L_{\infty})(\mathcal{M})$ if and only
if $\mu(A)\in (L_{1}\cap L_{\infty})(\mathbb{R_+}).$ Moreover,
$$\|A\|_{(L_{1}\cap L_{\infty})(\mathcal{M})}=\|\mu(A)\|_{(L_{1}\cap L_{\infty})(\mathbb{R_+})}, \,\,\,\, A\in (L_{1}\cap L_{\infty})(\mathcal{M}). $$
For an in-depth description on these spaces, refer to \cite[Chapter III]{DdePS} (see also \cite{DDdeP1}).
As considered in the commutative case, every symmetric Banach operator space satisfies
\begin{equation}\label{imbed}
(L_1\cap L_{\infty})(\mathcal{M})\subset E(\mathcal{M})\subset (L_1+L_{\infty})(\mathcal{M})
\end{equation}
(see \cite[Theorem 4.5, pp. 36-37]{DdePS} for more details).

\subsection{Lorentz spaces}\label{NC Lorentz}

\begin{definition}\cite[Definition II. 1.1, p.49]{KPS} A function $\varphi$ on the semiaxis $\mathbb{R}_{+}$ is said to be quasiconcave if
\begin{enumerate}[{\rm (i)}]
         \item $\varphi(t)=0\Leftrightarrow t=0;$
          \item $\varphi(t)$ is positive and increasing for $t>0;$
         \item $\frac{\varphi(t)}{t}$ is decreasing for $t>0.$
\end{enumerate}
\end{definition}
Observe that every nonnegative function on $[0,\infty)$ that vanishes only at origin is quasiconcave. The reverse, however, is not always true.
However, we may replace, if necessary, a quasiconcave function $\varphi$ by its least concave majorant $\widetilde{\varphi}$ such that
$$\frac{1}{2}\widetilde{\varphi}\leq\varphi\leq\widetilde{\varphi}$$
(see \cite[Proposition 5.10, p. 71]{BSh}).

Let $\Omega$ denote the set of increasing concave functions $\varphi:[0,\infty)\rightarrow[0,\infty)$ for which $\lim_{t\rightarrow 0+}\varphi(t)=0$ (or simply $\varphi(+0)=0$).
For the function $\varphi$ in $\Omega,$ the Lorentz space $\Lambda_{\varphi}(\mathbb{R_+})$ is defined by setting
$$\Lambda_{\varphi}(\mathbb{R_+}):=\left\{x\in S(\mathbb{R_+}): \int_{\mathbb{R_+}}\mu(s,x)d\varphi(s)<\infty\right\},$$
and equipped with the norm
$$\|x\|_{\Lambda_\varphi(\mathbb{R_+})}:=\int_{\mathbb{R_+}}\mu(s,x)d\varphi(s).$$
These spaces are examples of symmetric Banach function spaces.
For more details on Lorentz spaces, we refer the reader to \cite[Chapter II.5]{BSh} and \cite[Chapter II.5]{KPS}.

As in the commutative case, for a function $\varphi$ in $\Omega$ define the corresponding non-commutative Lorentz space by setting
$$\Lambda_{\varphi}(\mathcal{M}):=\left\{A\in S(\mathcal{M},\tau): \int_{\mathbb{R_+}}\mu(s,A)d\varphi(s)<\infty\right\}$$
equipped with the norm
\begin{equation}\label{Lphi}
\|A\|_{\Lambda_\varphi(\mathcal{M})}:=\int_{\mathbb{R_+}}\mu(s,A)d\varphi(s).
\end{equation}
These operator spaces become symmetric operator spaces.

Let $\psi$ be a quasiconcave function on $(0,\infty)$. The space
$$M_{\psi}(\mathbb{R_+})=\{f \in S(\mathbb{R_+}): \|f\|_{M_{\psi}} < \infty \}$$
equipped with the norm
$$\|f\|_{M_{\psi}(\mathbb{R_+})}=\sup_{t>0} \frac{\psi(t)}{t} \cdot \int_{0}^{t} \mu(s,f)ds$$
is the symmetric space with the fundamental function $\varphi(t)=\frac{t}{\varphi(t)}\cdot \chi_{(0,\infty)}(t).$ The space $(M_{\psi}, \| \cdot\|_{M_{\psi}})$ is called the Marcinkiewicz space.

Similarly to the Lorentz space, define the non-commutative Marcinkiewicz space
$$M_{\psi}(\mathcal{M}):=\left\{A\in S(\mathcal{M},\tau): \sup_{t>0} \frac{\psi(t)}{t} \cdot \int_{0}^{t} \mu(s,A)ds<\infty\right\}$$
equipped with the norm
$$\|A\|_{M_{\psi}(\mathcal{M})}=\sup_{t>0} \frac{\psi(t)}{t} \cdot \int_{0}^{t} \mu(s,A)ds.$$

\section{Upper estimate of generalized singular number of the non-commutative Hardy-Littlewood maximal function}
In this section we estimate generalised singular number of the non-commutative Hardy-Littlewood maximal function. We also show interesting applications of this result.

To prove the principal result the following results are necessary.
\begin{lem}\label{tau comp} Every operator $A$ in $\mathcal{L}_{loc}(\mathcal{M},\tau)$ is $\tau$-compact.
\end{lem}
\begin{proof} Assume, without loss of generality, that $0\leq A\in \mathcal{L}_{loc}(\mathcal{M},\tau).$ Since $A$ is $\tau$-measurable, then $\tau(E_{(N,\infty)}(A))<\infty$. Fix $\varepsilon>0$ and consider
$\tau(E_{(\varepsilon,N)}(A)).$ Since $\frac{1}{\varepsilon}AE_{(\varepsilon,N)}(A)>1,$ it follows from \eqref{*} that
$$\tau(E_{(\varepsilon,N)}(A))\leq \frac{1}{\varepsilon}\tau(AE_{(\varepsilon,N)}(A))<\infty.$$
Thus, $\tau(E_{(\varepsilon,\infty)}(A))<\infty$ for any $\varepsilon>0.$ Hence, $A$ is $\tau$-compact.
\end{proof}

Let $C:L_{1}(\mathbb{R}_+)\to L_{1,\infty}(\mathbb{R}_+)$ be the Ces\`{a}ro operator defined by

\begin{equation}\label{cesaro}
(Cf)(t):=\frac{1}{t}\int_{0}^{t}f(s)ds, \quad f\in L_{1}(\mathbb{R}_+), \quad t>0.
\end{equation}

The following theorem is the main result of this section.
\begin{thm}\label{main th}
For every $A\in \mathcal{L}_{loc}(\mathcal{M},\tau),$ we have
\begin{equation}\label{16}
	\mu(t, MA(|A|))\leq 16\cdot (C\mu(A))(t), \ \ \ \forall t>0.
\end{equation}	
\end{thm}

\begin{proof} The idea of the proof is similar to the commutative case.
	
	Let $A\in \mathcal{L}_{loc}(\mathcal{M},\tau).$ Then by Lemma \ref{tau comp}, $A$ is a $\tau$-compact operator. First, we prove the following inequality
\begin{equation}\label{equal}
	\|A_{1}\|_{\mathcal{L}_{1}(\mathcal{M},\tau)}+t\|A_{2}\|_{\mathcal{L}_{\infty}(\mathcal{M},\tau)}\leq 2t(C\mu(A))(t), \ \ \ t > 0. \end{equation}
 Fix $\xi_0>0$, suppose $(C\mu(A))(t)<+\infty,$ otherwise there is nothing to prove.
 Let $\xi_0=\mu(\xi,A).$
Set $A_{1}$ and $A_{2}$ as follows $A_1=A \cdot \chi_{(0, \xi_0]}$ and $A_2=A-A_1=A \cdot \chi_{(\xi_0,\infty)}$.
Then $||A_1||_{L_1(\mathcal{M},\tau)}=\int_{0}^{\xi_0}\mu(s,A)ds$ and
$$||A_2||_{L_\infty(\mathcal{M},\tau)}=\text{ess}\sup_{s\geq 0} |\mu(s+\xi_0),A|\leq \mu(\xi_0,A).$$
Since $\mu(s,A)$ is decreasing, we obtain
$$\xi_0\cdot ||A_2||_{L_{\infty}(\mathcal{M},\tau)}\leq \xi_0 \cdot \mu(\xi_0,A)\leq \int_{0}^{\xi_0}\mu(s,A)ds.$$
Therefore,
$$\|A_{1}\|_{\mathcal{L}_{1}(\mathcal{M},\tau)}+\xi_0\|A_{2}\|_{\mathcal{L}_{\infty}(\mathcal{M},\tau)}\leq 2\int_{0}^{\xi_0}\mu(s,A)ds=2\xi_0(C\mu(A))(\xi_0).$$ Since $\xi_0$ is arbitrary, the inequality (\ref{equal}) holds.




Let $A_1$ and $A_2$ be as defined above. Then, it was shown in the proof of the Lemma 3.1 in \cite{Sh} that
\begin{eqnarray*}
  \frac{1}{\tau(E_{[x-r,x+r]}(|A|))}\tau(|A|E_{[x-r,x+r]}(|A|)) &\leq & \frac{1}{\tau(E_{[x-r,x+r]}(|A|))}\tau(|A_1|E_{[x-r,x+r]}(|A_1|))\\
  &+& \frac{1}{\tau(E_{[x-r,x+r]}(|A|))}\tau(|A_2|E_{[x-r,x+r]}(|A_2|))\\
&\leq&\frac{1}{\tau(E_{[x-r,x+r]}(|A_1|))}\tau(|A_1|E_{[x-r,x+r]}(|A_1|))\\
  &+& \frac{1}{\tau(E_{[x-r,x+r]}(|A_2|))}\tau(|A_2|E_{[x-r,x+r]}(|A_2|)).
\end{eqnarray*}
Hence, taking supremum over $r>0,$ we obtain
$$MA(x)\leq MA_1(x)+MA_2(x),$$
which implies that
\begin{equation}\label{sublinear}
  MA(|A|)\leqslant MA_1(|A_1|)+MA_2(|A_2|).
\end{equation}
Then, by \eqref{sublinear} and  Lemma 2.5 in \cite{FK} we have
\begin{eqnarray*}
  \mu(s,MA(|A|)) &\stackrel{\eqref{sublinear}}{\leq }& \mu(s,MA_1(|A_1|)+MA_2(|A_2|)) \\
  &\leq& \mu(\frac{s}{2},MA_{1}(|A_1|))+\mu(\frac{s}{2},MA_{2}(|A_2|)) \quad s>0.
\end{eqnarray*}
We estimate each term separately. For the first term, since $A_{1}\in \mathcal{L}_{1}(\mathcal{M},\tau),$ it follows from \cite[Theorem 1]{Bek} that
\begin{equation}\label{L1}\mu(\frac{s}{2},MA_{1}(|A_1|))\leq \frac{16}{s}\|A_1\|_{\mathcal{L}_{1}(\mathcal{M},\tau)},
\end{equation}
and since $A_{2}\in \mathcal{L}_{\infty}(\mathcal{M},\tau),$ by
\cite[Lemma 1 (ii)]{Bek} we have
\begin{equation}\label{L}\mu(\frac{s}{2},MA_{2}(|A_2|)\leq \mu(0,MA_{2}(|A_2|)=\|MA_{2}(|A_2|)\|_{\mathcal{L}_{\infty}(\mathcal{M},\tau)}\leq \|A_2\|_{\mathcal{L}_{\infty}(\mathcal{M},\tau)}.
\end{equation}
Combining \eqref{L1} and \eqref{L} we infer
\begin{eqnarray*}
  \mu(s,MA(|A|)) \leq \frac{16}{s}\|A_1\|_{\mathcal{L}_{1}(\mathcal{M},\tau)}+\|A_2\|_{\mathcal{L}_{\infty}(\mathcal{M},\tau)} \quad s>0.
\end{eqnarray*}

Putting $s=t$ and using \eqref{equal}, we obtain
\begin{eqnarray*}\mu(t,MA(|A|))&\leq&\frac{16}{t}\|A_1\|_{\mathcal{L}_{1}(\mathcal{M},\tau)}+\|A_2\|_{\mathcal{L}_{\infty}(\mathcal{M},\tau)}\\
&\leq&\frac{16}{t}\|A_1\|_{\mathcal{L}_{1}(\mathcal{M},\tau)}+16\|A_2\|_{\mathcal{L}_{\infty}(\mathcal{M},\tau)}\\
&\leq& \frac{16}{t}\left(\|A_{1}\|_{\mathcal{L}_{1}(\mathcal{M},\tau)}+t\|A_{2}\|_{\mathcal{L}_{\infty}(\mathcal{M},\tau)}\right)\stackrel{\eqref{equal}}{\leq}16\cdot(C\mu(A))(t).
\end{eqnarray*}
Since $A$ is arbitrary, we obtain desired inequality.\end{proof}

As a corollary we obtain the results of Lemma 1 (ii) in \cite{Bek}.
\begin{corollary}\label{infty} For any $A\in \mathcal{L}_{\infty}(\mathcal{M},\tau),$ we have
$$\|MA(|A|)\|_{\mathcal{L}_{\infty}(\mathcal{M},\tau)}\leq 16\cdot\|A\|_{\mathcal{L}_{\infty}(\mathcal{M},\tau)}.$$
\end{corollary}
\begin{proof}
It is well-known that $C$ is bounded operator from $L_{\infty}(\mathbb{R}_{+})$ into itself, that is,
\begin{equation}\label{Hardy oper}\|Cf\|_{L_{\infty}(\mathbb{R}_{+})}\leq \|f\|_{L_{\infty}(\mathbb{R}_{+})}, \quad \forall f\in {L_{\infty}(\mathbb{R}_{+})}.
\end{equation}
Therefore, it follows from Theorem \ref{main th} that
\begin{eqnarray*}\|MA(|A|)\|_{\mathcal{L}_{\infty}(\mathcal{M},\tau)}&=&\|\mu(MA(|A|))\|_{L_{\infty}(\mathbb{R}_{+})}\\
&\leq&16\cdot\|C\mu(A)\|_{L_{\infty}(\mathbb{R}_{+})}\\
&\stackrel{\eqref{Hardy oper}}{\leq}&16\cdot\|\mu(A)\|_{L_{\infty}(\mathbb{R}_{+})}=16\cdot\|A\|_{\mathcal{L}_{\infty}(\mathcal{M},\tau)}, \ \ \forall A \in L_{\infty}(\mathcal{M}, \tau).
\end{eqnarray*}
\end{proof}

Similarly, we complement the main results of \cite{Sh} (see \cite[Theorem 3.2]{Sh}) for the general case when $\mathcal{M}$ is a semi-finite von Neumann with a normal faithful semi-finite trace $\tau$.
\begin{corollary}\label{lorentz type} Let $0< q<\infty$, $1<p,p_0,p_1<\infty$ and $p_0\neq p_1$ be such that
$\frac{1}{p}=\frac{1-\theta}{p_0}+\frac{\theta}{p_1}$
for some $0<\theta < 1.$
Then there exists a constant $c_{pq}>0$ such that for all
$A\in \mathcal{L}_{p,q}(\mathcal{M},\tau)$ we have
$$\|MA(|A|)\|_{\mathcal{L}_{p,q}(\mathcal{M},\tau)}\leq c_{p,q}\cdot\|A\|_{\mathcal{L}_{p,q}(\mathcal{M},\tau)}.$$
\end{corollary}
\begin{proof}By Hardy's inequality in \cite{Hunt} (see also \cite[Lemma 3.9, p.124]{BSh}), we have
\begin{equation}\label{Hardy ineq}\|Cf\|_{L_{p,q}(\mathbb{R}_{+})}\leq \widetilde{c}_{p,q}\cdot\|f\|_{L_{p,q}(\mathbb{R}_{+})}, \quad \forall f\in {L_{p,q}(\mathbb{R}_{+})}.
\end{equation}
By Theorem \ref{main th}, we obtain
\begin{eqnarray*}\|MA(|A|)\|_{\mathcal{L}_{p,q}(\mathcal{M},\tau)}&=&\|\mu(MA(|A|))\|_{L_{p,q}(\mathbb{R}_{+})}\\
&\leq&16\cdot\|C\mu(A)\|_{L_{p,q}(\mathbb{R}_{+})}\\
&\stackrel{\eqref{Hardy ineq}}{\leq}&16\cdot \widetilde{c}_{p,q}\cdot\|\mu(A)\|_{L_{p,q}(\mathbb{R}_{+})}=c_{p,q}\cdot\|A\|_{\mathcal{L}_{p,q}(\mathcal{M},\tau)},
\end{eqnarray*}
where $c_{p,q}=8\cdot\widetilde{c}_{p,q}.$
This concludes the proof.
\end{proof}
In particular, when $p=q,$ we obtain the result of Theorem 2 in \cite{Bek}.
\begin{corollary}\label{strong type} For $1<p<\infty$ and $A\in \mathcal{L}_p(\mathcal{M},\tau),$ there is a constant $c_p> 0$ such that
$$\|MA(|A|)\|_{\mathcal{L}_{p}(\mathcal{M},\tau)}\leq c_p\cdot\|A\|_{\mathcal{L}_{p}(\mathcal{M},\tau)}.$$
\end{corollary}
\begin{proof}The proof follows from the previous corollary.
\end{proof}
Let us denote
$$(L_{1,\infty}(\mathbb{R}_{+}))^{0}:=\{x\in L_{1,\infty}(\mathbb{R}_{+}):\lim_{t\rightarrow0+} t\mu(t,x)=0\}.$$
Note that this space coincides with the closure of all bounded functions in $L_{1,\infty}(\mathbb{R}_+).$

The following result refines the result of Theorem 1 in \cite{Bek}.
\begin{thm}\label{special case} The non-commutative Hardy-Littlewood maximal function
$$MA(|\cdot|):\mathcal{L}_{1}(\mathcal{M},\tau) \to (\mathcal{L}_{1,\infty}(\mathcal{M},\tau))^{0}$$ is bounded.	
\end{thm}
\begin{proof} It was shown in \cite[Remark 3.5.]{ST2} that $C:L_{1}(\mathbb{R}_+)\to (L_{1,\infty}(\mathbb{R}_+))^{0}$ is bounded.
 Hence, by Theorem \ref{main th} we obtain the desired result. Indeed, let $A\in \mathcal{L}_{1}(\mathcal{M},\tau).$ Then by the definition $\mu(A)\in L_{1}(\mathbb{R}_+)$
Since $C\mu(A)\in(L_{1,\infty}(\mathbb{R}_+))^{0},$ it follows from Theorem \ref{main th} that $\mu(MA(|A|))\in (L_{1,\infty}(\mathbb{R}_+))^{0},$ that means
$MA(|A|)\in (\mathcal{L}_{1,\infty}(\mathcal{M},\tau))^{0}.$
Since $A$ is arbitrary, this completes the proof.
\end{proof}

\section{The Hardy-Littlewood maximal function on non-commutative symmetric spaces}
Let $E(\mathbb{R_+})$ be a symmetric space of functions.
In \cite[Theorem 2.3]{ST2}, \cite[Theorem 3.5]{ST1} the authors provide the minimal symmetric space $F(\mathbb{R_+})$ with Fatou norm such that the Ces\`{a}ro operator (defined as in (\ref{cesaro}))
$C: E(\mathbb{R_+}) \to F(\mathbb{R_+})$ is bounded.
Following their terminology, denote by $\mathcal{R}[C, E]$ the minimal symmetric space $F(\mathbb{R_+})$ such that
$$C: E(\mathbb{R_+}) \to F(\mathbb{R_+})$$ is bounded.
Let $\varphi$ be a quasiconcave function satisfying
\begin{equation}\label{varphi condition}
\varphi(t) \geq c t\log(1+1/t),\ c>0.
\end{equation}
Then the Ces\`{a}ro operator (see \cite[Theorem 2.3]{ST2})
$$C: \Lambda_{\varphi}(\mathbb{R_+}) \to F(\Lambda_{\varphi}):=\left\{f\in (L_1+L_{\infty})(\mathbb{R_+}): \exists g\in\Lambda_{\varphi}(\mathbb{R_+}) \ \ \text{such that} \ \ \mu(f)\prec\prec C\mu(g)\right\}$$
is bounded.

Define
\begin{equation}\label{F}
F(\mathbb{R_+}):=\{f\in (L_1+L_{\infty})(\mathbb{R_+}): \exists g\in E(\mathbb{R_+}) \ \ \text{such that} \ \ \mu(f) \prec\prec C\mu(g)\}.
\end{equation}
It is known that the Ces\`{a}ro operator
$$C: E(\mathbb{R_+}) \to F(\mathbb{R_+})$$
is bounded. Moreover, such defined $F(\mathbb{R_+})$ is minimal among all symmetric spaces.
Now we define non-commutative symmetric spaces. It was proved  (see \cite[Theorem 3.1.1]{LSZ}, see also \cite{KS}) that if $E$ is a symmetric Banach function space on $\mathbb{R_+}$, then the corresponding non-commutative symmetric space of $\tau-$measurable operators, defined as
$$E(\mathcal{M}):=\{A\in S(\mathcal{M}): \mu(A) \in E(\mathbb{R_+})\},$$
is also a symmetric space equipped with the norm
$$||A||_{E(\mathcal{M})}:=||\mu(A)||_{E(\mathbb{R_+})}.$$

Similarly, we define
$$F(\mathcal{M}):=\{A\in S(\mathcal{M}): \mu(A) \in F(\mathbb{R_+})\},$$
equipped with the norm
$$||A||_{F(\mathcal{M})}:=||\mu(A)||_{F(\mathbb{R_+})},$$
where $F$ is defined as in \eqref{F}.

The following theorem is the main result of this section.
\begin{thm}\label{main1}
The non-commutative Hardy-Littlewood maximal function
$$MA(\cdot): E(\mathcal{M}) \to F(\mathcal{M})$$
is bounded.
\end{thm}
\begin{proof}

Take $A \in E(\mathcal{M})$. Then $\mu(A) \in E(\mathbb{R_+})$. Since $C: E(\mathbb{R_+}) \to F(\mathbb{R_+})$ is bounded, it follows from the Theorem 1.3 that
\begin{eqnarray*}
||MA(|A|)||_{F(\mathcal{M})}=||\mu(MA(|A|))||_{F(\mathbb{R_+})}
{\leq}&16||C\mu(A)||_{F(\mathbb{R_+})}
\lesssim  ||\mu(A)||_{E(\mathbb{R_+})}=||A||_{E(\mathcal{M})}
\end{eqnarray*}
Since $A$ is arbitrary the assertion follows.
\end{proof}

As an immediate corollary we obtain:
\begin{proposition}\label{16}
Let $\varphi$ be a quasi-concave function satisfying (\ref{varphi condition}) and let $\psi$ be a quasi-concave function such that
\begin{equation}\label{psi varphi}
\int_{t}^{\infty}\frac{\psi(s)}{s^2}ds \lesssim \frac{\varphi(t)}{t}.
\end{equation}
Then the non-commutative Hardy-Littlewood maximal function
$$MA(\cdot): \Lambda_{\varphi}(\mathcal{M}) \to \Lambda_{\psi}(\mathcal{M})$$
is bounded.
\end{proposition}
\begin{proof}
It is known that if $\varphi$ satisfies (\ref{varphi condition}), the Ces\`{a}ro operator (see \cite[Proposition 4.4]{ST1})
$$C: \Lambda_{\varphi}(\mathbb{R_+}) \to \Lambda_{\psi}(\mathbb{R_+})$$
is bounded if and only if
(\ref{psi varphi}) holds.
Moreover, $\Lambda_{\psi}(\mathbb{R_+})$ is minimal among such symmetric Banach function spaces. Therefore, similarly to the proof of the Theorem \ref{main1}, it follows from the Theorem \ref{main th} that the non-commutative Hardy-Littlewood maximal function
$$MA(\cdot): \Lambda_{\varphi}(\mathcal{M}) \to \Lambda_{\psi}(\mathcal{M})$$
is bounded.
\end{proof}

Now we consider the boundedness of the Hardy-Littlewood maximal function on non-commutative Marcinkiewicz spaces, which are the dual spaces to the Lorentz spaces.

Similarly to the Lorentz space, we obtain the following boundedness criterion for the Marcinkiewicz space.
\begin{proposition}\label{17}
Let $\phi$ be a quasi-concave function such that $1/\phi$ is locally integrable at zero.
Then the non-commutative Hardy-Littlewood maximal function
$$MA(\cdot): M_{\phi}(\mathcal{M}) \to M_{\psi}(\mathcal{M})$$
is bounded, where
\begin{equation}\label{inv}
\psi(t)=t\cdot \left( \int_{0}^{t}\frac{ds}{\phi(s)}\right)^{-1}, \ t>0.
\end{equation}
\end{proposition}

\begin{proof}
In \cite[Theorem 4.7]{ST1} it was proved that the Ces\`{a}ro operator
$C: M_{\phi}(\mathbb{R_+}) \to M_{\psi}(\mathbb{R_+})$ is bounded if and only if $1/\phi$ is locally integrable at zero. Moreover, in this case the optimal range is $M_{\psi}(\mathbb{R_+}).$
\end{proof}

$\mathbf{Example \ 1.}$
Let $\varphi(t)=t \log^2(1+\frac{1}{\sqrt{t}}).$ Then it is easy to show that $\varphi$ satisfies (\ref{varphi condition}).
Take $\phi(t)=t \log(1+\frac{1}{t}).$ Then
$$\int_{t}^{\infty}\frac{\phi(s)}{s^2}ds \lesssim \frac{\varphi(t)}{t}, \ \ t>0.$$

 Indeed, for $0<t \leq 1:$
 \begin{eqnarray*}
 \begin{split}
  \int_{t}^{\infty}\frac{\phi(s)}{s^2}ds&=\int_{t}^{\infty}\frac{s\log(1+\frac{1}{s})}{s^2}ds=\int_{t}^{\infty}\frac{1}{s}\log(1+\frac{1}{s})ds\\ &\approx \int_{t}^{1}\frac{1}{s}\log(1+\frac{1}{s})ds +1 \leq \log(1+\frac{1}{t})\log(\frac{1}{t})+1\\
  &\leq c_{abs} \log^2(1+\frac{1}{t})=c_{abs}\log^2(1+\frac{1}{t})^{1/2}\\
  &\leq c_{abs}\log^2(1+\frac{1}{\sqrt{t}})=c_{abs} \cdot \frac{t \log^2(1+\frac{1}{\sqrt{t}})}{t}=c_{abs} \cdot \frac{\varphi(t)}{t}.
 \end{split}
 \end{eqnarray*}

 For $t>1:$
 $$\int_{t}^{\infty}\frac{\phi(s)}{s^2}ds=\int_{t}^{\infty}\frac{s\log(1+\frac{1}{s})}{s^2}ds=\int_{t}^{\infty}\frac{1}{s}\log(1+\frac{1}{s})ds$$
 $$\approx\int_{t}^{\infty} \frac{ds}{s^2}=\frac{1}{t}\approx \log^2(1+\frac{1}{\sqrt{t}})=c_{abs} \cdot \frac{t \log(1+\frac{1}{\sqrt{t}})}{t}=c_{abs}\cdot \frac{\varphi(t)}{t}.$$
 Hence, by Proposition \ref{16} we obtain that
$$MA(\cdot): \Lambda_{t\log^2(1+\frac{1}{\sqrt{t}})}(\mathcal{M}) \to \Lambda_{t\log(1+\frac{1}{t})}(\mathcal{M})$$
is bounded.

$\mathbf{Example \ 2.}$
Let us consider the function $\phi(t)=max\{1,t\}, t>0.$ Then by (\ref{inv}), we obtain
$\psi(t)\approx\frac{t}{\log(1+t)}.$ Since $M_{\phi}(\mathcal{M})=(L_1+L_{\infty})(\mathcal{M})$ for the function $\phi(t)=max\{1,t\}$, it follows from Proposition \ref{17} that
$$MA(\cdot): (L_1+L_\infty)(\mathcal{M}) \to M_{\frac{t}{\log(1+t)}}(\mathcal{M})$$
is bounded.

\section{Acknowledgment}
Authors would like to thank Professor F. Sukochev for his very useful comments and helps to prove Lemma \ref{tau comp}.
 The work was partially supported by the grants (No. AP08052004 and No. AP08051978) of the Science Committee of the Ministry of Education and Science of the Republic of Kazakhstan.


\begin{thebibliography}{99}

\bibitem{Bek} T.N. Bekjan,  {\it Hardy-Littlewood maximal function of $\tau$-measurable operators}, J. Math. Anal. Appl. {\bf 322} (2006), 87--96.

\bibitem{BSh} C. Bennett and R. Sharpley, {\it Interpolation of Operators}, Pure and Applied Mathematics, {\bf 129}. Academic Press, 1988.


\bibitem{DDdeP1}  P.G. Dodds, T.K. Dodds, B. de Pagter, {\it Noncommutative K\"{o}the duality},  Trans. Amer. Math. Soc., {\bf 339} (1993), 717--750.

\bibitem{DdePS} P.G. Dodds, B. de Pagter, and F.A. Sukochev, {\it Theory of Noncommutative integration}, Unpublished manuscript.

\bibitem{FK} T. Fack, H. Kosaki, {\it Generalized $s$-numbers of $\tau$-measurable operators}, Pacific J. Math., {\bf 123(2)} (1986), 269--300.

\bibitem{Gra} L.Grafakos, {\it Classical Fourier analysis}, 2nd ed., Springer, Berlin, (2008).


\bibitem{Junge} M.Junge, {\it Doob's inequality for non-commutative martingales}, J.Reine Angew.Math. 549, (2002), 149-190.

\bibitem{Xu} M.Junge, Q.Xu, {\it Non-commutative maximal ergodic theorems}, J.Am.Math.Soc. 20(2), (2007), 385-439.

\bibitem{Hunt} R.A. Hunt, {\it On $L(p,q)$ spaces}. Enseign. Math. {\bf 12} (1966), 249--276.


\bibitem{KS} N. Kalton, F. Sukochev, {\it Symmetric norms and spaces of operators.} J. Reine Angew. Math. {\bf 621} (2008), 81--121.


\bibitem{KPS} S. Krein, Y. Petunin, and E. Semenov, {\it Interpolation of linear operators}, Amer. Math. Soc., Providence, R.I., (1982).

\bibitem{LT} J. Lindenstrauss, L. Tzafiri, {\it Classical Banach spaces}. Springer-Verlag, II, (1979).

\bibitem{LSZ} S. Lord, F. Sukochev, D. Zanin, {\it Singular traces. Theory and applications}, De Gruyter Studies in Mathematics, {\bf 46}. De Gruyter, Berlin, 2013.

\bibitem{Mei} T.Mei, {\it Operator valued Hardy spaces}, Mem.Am.Math.Soc. 188, 881, (2007).


\bibitem{Sh} J. Shao, {\it Hardy-Littlewood maximal function on noncommutative Lorentz spaces}, Journal of Inequalities and Applications 2013, 2013:384.

\bibitem{ST2} J. Soria, P. Tradacete, {\it Characterization of the restricted type spaces $R(X)$}, Math. Ineq. Applic. {\bf 18} (2015), 295--319.

\bibitem{ST1} J. Soria, P. Tradacete, {\it Optimal rearrangement invariant range for Hardy-type operators }, Proc. of the Royal Soc. of Edinburgh, {\bf 146A} (2016), 865--893.

\bibitem{SW} E.Stein, G.Weiss, {\it Introdunction to Fourier analysis on Euclidean spaces}, Princeton University Press, Princeton (1971).

\bibitem{F} F. Sukochev, {\it Completeness of quasi-normed symmetric operator spaces.} Indag. Math. (N.S.) {\bf 25} (2014), no. 2, 376--388.

\bibitem{STZ1} F. Sukochev, K. Tulenov, and D. Zanin, {\it The optimal range of the Calder\'{o}n operator and its applications}, J. Func. Anal. {\bf 277:10} (2019), 3513--3559.

\bibitem{STZ2} F. Sukochev, K. Tulenov, and D. Zanin, {\it Nehari-Type Theorem for Non-commutative Hardy
Spaces}, J. Geom. Anal. {\bf 77} (2017), 1789--1802.

\end{thebibliography}
\end{document}